\documentclass[a4paper,conference]{IEEEtran}
\usepackage{amsmath}
\usepackage{amsthm}
\usepackage{amssymb}
\usepackage{bbm}
\usepackage{verbatim}
\usepackage{graphicx}
\usepackage{subfigure}
\usepackage{afterpage}
\usepackage{etoolbox}
\usepackage{footbib}
\usepackage{float}
\usepackage{algorithmic}
\usepackage[algosection,algoruled,vlined]{algorithm2e}
\usepackage{cite}
\usepackage{url}

\numberwithin{equation}{section}
\numberwithin{figure}{section}
\numberwithin{table}{section}

\makeatletter
\renewcommand{\p@subfigure}{\thefigure}
\makeatother

\newtheorem{theorem}{Theorem}[section]

\newtheorem{lemma}[theorem]{Lemma}

\newcommand{\repeatable}[2]{\makeatletter \global\expandafter\def\csname repText@#1\endcsname {#2} \makeatother #2}
\newcommand{\repeatxt}[1]{\makeatletter \expandafter\csname repText@#1\endcsname \makeatother}

\newtoggle{citeInAbstract}
\toggletrue{citeInAbstract}

\newcommand{\usecrop}[2]
{
	\newlength{\cropwidth}
	\setlength{\cropwidth}{\the\textwidth}
	\addtolength{\cropwidth}{#1}
	\newlength{\cropheight}
	\setlength{\cropheight}{\the\textheight}
	\addtolength{\cropheight}{#2}
	\usepackage[width=\the\cropwidth,height=\the\cropheight,center]{crop}
}

\DeclareMathAlphabet{\mathpzc}{OT1}{pzc}{m}{it}

\begin{document}

\title{Missing Entries Matrix Approximation and Completion}

\author{\IEEEauthorblockN{Gil Shabat}
\IEEEauthorblockA{School of Electrical Engineering\\
Tel Aviv University\\
gil@eng.tau.ac.il}
\and
\IEEEauthorblockN{Yaniv Shmueli}
\IEEEauthorblockA{School of Computer Science\\
Tel Aviv University\\
yaniv.shmueli@cs.tau.ac.il}
\and
\IEEEauthorblockN{Amir Averbuch}
\IEEEauthorblockA{School of Computer Science\\
Tel Aviv University\\
amir@math.tau.ac.il}
}

\maketitle

\begin{abstract}
We describe several algorithms for matrix completion and matrix approximation when only some of its entries are known. 
The approximation constraint can be any whose approximated solution is known for the full matrix. For low rank approximations, similar algorithms appears recently in the literature under different names. In this work, we introduce new theorems for matrix approximation and show that these algorithms can be extended to handle different constraints such as nuclear norm, spectral norm, orthogonality constraints and more that are different than low rank approximations. As the algorithms can be viewed from an optimization point of view, we discuss their convergence to global solution for the convex case. We also discuss the optimal step size and show that it is fixed in each iteration. In addition, the derived matrix completion flow is robust and does not require any parameters. This matrix completion flow is applicable to different spectral minimizations and can be applied  to physics, mathematics and electrical engineering problems such as data reconstruction of images and data coming from PDEs such as Helmholtz's equation used for electromagnetic waves.
\end{abstract}

\section{Introduction}
Matrix completion and matrix approximation are important problems in
a variety of fields such as statistics \cite{sta}, biology
\cite{ml}, statistical machine learning \cite{SJ}, signal processing and
computer vision/image processing \cite{mikielad}. Rank
reduction by matrix approximation is important, for example, in
compression where low rank indicates the existence of redundant
information and matrix completion is important in collaborative filtering,
such as the Netflix problem and different reconstruction problems.
Usually, the matrix completion problem, is defined as finding a matrix, with 
smallest possible rank, that satisfy the existence of certain entries.
\begin{equation}
\label{lowrank_eq}
\begin{array}{l}
 \mbox{minimize   } \mbox{rank }(\mathbf{X}) \\
\mbox{subject to } X_{i,j}=M_{i,j},\;\;\; (i,j) \in \Omega .
\end{array}
\end{equation}
Since Eq.~\ref{lowrank_eq} is an NP-hard problem, some
relaxations methods have been proposed. The most popular relaxation
is one that replaces the rank by the nuclear norm:
\begin{equation}
\begin{array}{l}
\mbox{minimize   } \Vert\mathbf{X}\Vert_* \\
\mbox{subject to } X_{i,j}=M_{i,j},\;\;\; (i,j) \in \Omega ,
\end{array}
\end{equation}
where $\Vert\mathbf{X}\Vert_*$ denotes the nuclear norm of $\mathbf{X}$
that is equal to the sum of the singular values of $\mathbf{X}$. A
small value of $\Vert\mathbf{X}\Vert_*$ is related to the property of
having a low rank \cite{fazel}. An iterative solution, which is
based on a singular value thresholding, is given in \cite{candes}. A
completion algorithm, based on the local information of the matrix,
is proposed in \cite{Nan}. In this work, a more robust and simple approach for solving
a variety of matrix approximation of certain entries by approximating the full matrix
is discussed. We approximate problems of the form
\begin{equation}
\label{eq1}
\begin{array}{l}
\mbox{minimize   } \Vert \mathcal{P}_{\Omega}\mathbf{X}-\mathcal{P}_{\Omega}\mathbf{M} \Vert_F \\
\mbox{subject to } f(\mathbf{X}) \le 0,
 \end{array}
\end{equation}
given that the solution for
\begin{equation}
\label{eq2}
\begin{array}{l}
\mbox{minimize   } \Vert \mathbf{X}-\mathbf{M} \Vert_F \\
\mbox{subject to } f(\mathbf{X}) \le 0
 \end{array}
\end{equation}
is known. Here, $\{\mathcal{P}_{\Omega}\mathbf{X}\}_{i,j}=X_{i,j}$ if $(i,j) \in \Omega$ and $0$ otherwise.
If $f(\mathbf{X})$ is convex and satisfies some condition (which is explained in the next sections), the algorithm finds the global solution.
Nevertheless, convergence is guaranteed, but to a local solution.
Then, we show how this algorithm can be used for solving a 
variety of matrix completion problems as well, such as spectral norm completion:
\begin{equation}
\begin{array}{l}
\mbox{minimize   } \Vert\mathbf{X}\Vert_2 \\
\mbox{subject to } X_{i,j}=M_{i,j},\;\;\; (i,j) \in \Omega,
\end{array}
\end{equation}
Ky-Fan norm completion:
\begin{equation}
\begin{array}{l}
\mbox{minimize   } \Vert\mathbf{X}\Vert_{(k)} \\
\mbox{subject to } X_{i,j}=M_{i,j},\;\;\; (i,j) \in \Omega,
\end{array}
\end{equation}
where $\Vert\mathbf{X}\Vert_{(k)}=\sum_{i=1}^k \sigma_i$ (sum of largest $k$ singular values). 
Note that the spectral norm and the nuclear norm are a special case of the Ky-Fan norm.
We also discuss approximation problems such as:
\begin{equation}
\begin{array}{l}
\mbox{minimize   } \Vert\mathcal{P}_{\Omega}\mathbf{X}-\mathcal{P}_{\Omega}\mathbf{M}\Vert_F \\
\mbox{subject to } \mathbf{X^TX=I}.
\end{array}
\end{equation}

\section{Theorems on full matrix approximation}
\label{sec:related}
 The algorithm that approximates a matrix at
certain points requires from us to be able to approximate the matrix
when taking into account all its entries. Therefore, we review some
theorems on full matrix approximation theorems in addition to the
well known Eckart-Young theorem mentioned in the introduction. The
low rank approximation problem can be modified to approximate a
matrix under the Frobenius norm while having the Frobenius norm as a
constraint as well instead of having low rank. Formally,
\begin{equation}
\label{eq_ff}
\begin{array}{l}
\mbox{minimize   } \Vert \mathbf{X}-\mathbf{M} \Vert_F \\
\mbox{subject to } \Vert \mathbf{X} \Vert_F \le \lambda .
 \end{array}
\end{equation}
A solution for Eq. \ref{eq_ff} is given by
$\mathbf{X}=\mathbf{\frac{M}{\Vert M
\Vert_F}\min(\Vert\mathbf{M}\Vert_F, \lambda)}$.
\begin{IEEEproof}
The expression  $\Vert \mathbf{X} \Vert_F^2 \le \lambda^2$ can be
thought of as an $m \times n$ dimensional ball with radius $\lambda$
centered at the origin. $\mathbf{M}$ is an $m \times n$ dimensional
point. We are looking for a point $\mathbf{X}$ on the ball $\Vert
\mathbf{X} \Vert_F^2 = \lambda^2$ that has a minimal Euclidean
distance (Frobenius norm) from $\mathbf{M}$. If $\Vert \mathbf{M}
\Vert_F \le \lambda$ then  $\mathbf{X=M}$ and it is inside the ball
having a distance of zero. If $\Vert \mathbf{M} \Vert_F
> \lambda$, then the shortest distance is given by the line going
from the origin to $\mathbf{M}$ whose intersection with the sphere
$\Vert \mathbf{X} \Vert_F^2 \le \lambda^2$ is the closest point to
$\mathbf{M}$. This point is given by
$\mathbf{X}=\mathbf{\frac{M}{\Vert M \Vert_F}}\lambda$.
\end{IEEEproof}
An alternative approach uses the Lagrange multiplier in a
brute-force manner. This leads to a non-linear system of equations,
which are difficult to solve. Note that this problem can be easily
extended to the general case
\begin{equation}
\label{eq_pff}
\begin{array}{l}
\mbox{minimize   } \Vert \mathcal{P}\mathbf{X}-\mathcal{P}\mathbf{M} \Vert_F \\
\mbox{subject to } \Vert \mathbf{X} \Vert_F \le \lambda .
 \end{array}
\end{equation}

\begin{IEEEproof}
The proof is similar to the previous one but here we are looking for
a point $\mathbf{X}$ on the sphere that is the closest to a line
whose points $\mathbf{X'} \in \mathcal{H}$ satisfy
$\mathcal{P}\mathbf{X'}=\mathcal{P}\mathbf{M}$. By geometrical
considerations, this point is given by
$\mathbf{X}=\frac{\mathcal{P}\mathbf{M}}{\Vert \mathcal{P}\mathbf{M}
\Vert_F}\lambda$.
\end{IEEEproof}
\par\noindent
Hence, we  showed a closed form solution for the problem in Eq.
\ref{eq_pff}.
\par\noindent
Another example is the  solution to the problem:
\begin{equation}
\label{eq6}
\begin{array}{l}
\mbox{minimize   } \Vert \mathbf{X}-\mathbf{M} \Vert_F \\
\mbox{subject to } \mathbf{X}^T\mathbf{X}=\mathbf{I}.
 \end{array}
\end{equation}
This is known as the orthogonal Procrustes problem (\cite{peter})
and the solution is given by $\mathbf{X=UV^*}$, where the SVD of
$\mathbf{M}$ is given by $\mathbf{M}=\mathbf{U\Sigma V^*}$. The
solution can be extended to a matrix $\mathbf{X}$ satisfying
$\mathbf{X^TX=D^2}$, where $\mathbf{D}$ is a known or unknown
diagonal matrix. When $\mathbf{D}$ is unknown, the solution is the
best possible orthogonal matrix. When $\mathbf{D}$ is known, the
problem can be converted to become the orthonormal case (Eq.
\ref{eq6}) by substituting $\mathbf{X=VD}$ where $\mathbf{V^TV=I}$.
When $\mathbf{D}$ is unknown, the problem can be solved by applying
an iterative algorithm that is described in \cite{everson}.

We now examine the following problem:
\begin{equation}
\label{eq7}
\begin{array}{l}
\mbox{minimize   } \Vert \mathbf{X}-\mathbf{M} \Vert_F \\
\mbox{subject to } \Vert \mathbf{X} \Vert_2 \le \lambda .
 \end{array}
\end{equation}
A solution to this problem uses the Pinching theorem (\cite{bhatia}):

\begin{lemma}[Pinching theorem]
\label{pinching} For every matrix  $\mathbf{A}$ and a unitary matrix
$\mathbf{U}$ and for any norm satisfying $\Vert
\mathbf{UAU^*}\Vert=\Vert \mathbf{A}\Vert$ then $\Vert
\mbox{diag}(\mathbf{A}) \Vert \le \Vert \mathbf{A} \Vert$.
\end{lemma}
A proof is given in \cite{goldberg}. An alternative proof is given
in \cite{gil}.

\begin{lemma}[Minimization of the Frobenius norm under the spectral norm constraint]
Assume the SVD of $\mathbf{M}$ is given by $\mathbf{M=U\Sigma V^*}$
where $\mathbf{\Sigma}=\mbox{diag}(\sigma_1,..,\sigma_n)$. Then, the
matrix $\mathbf{X}$, which minimizes $\Vert \mathbf{X}-\mathbf{M}
\Vert_F$ such that $\Vert \mathbf{X} \Vert_2 \le \lambda$, is given
by $\mathbf{X=U\tilde{\Sigma}V^*}$ where $\tilde{\sigma_i}$ are the
singular values of $\tilde{\Sigma}$ and
$\tilde{\sigma_i}=\min(\sigma_i, \lambda), i=1, \ldots k, ~k
\le n$.
\label{spectral}
\end{lemma}

\begin{IEEEproof}
$\Vert \mathbf{X-M} \Vert_F = \Vert \mathbf{X-U\Sigma V^*}
\Vert_F=\Vert \mathbf{U^*XV-\Sigma} \Vert_F$. Since
$\mathbf{\Sigma}$ is diagonal,  $\Vert
\mathbf{\mbox{diag}(\mathbf{U^*XV})-\Sigma} \Vert_F \le \Vert
\mathbf{U^*XV-\Sigma} \Vert_F$. From Lemma \ref{pinching} we know
that $\Vert \mbox{diag}(\mathbf{U^*XV}) \Vert_2 \le \Vert
\mathbf{U^*XV} \Vert_2$. Therefore, $\mathbf{U^*XV}$ has to be
diagonal and the best minimizer under the spectral norm constraint
is achieved by minimizing each element separately yielding
$\mathbf{U^*XV}=\mbox{diag}(\min(\sigma_i,\lambda)),~i=1,\ldots
k, k \le n$. Hence, $\mathbf{X=U\tilde{\Sigma}V^*}$.
\end{IEEEproof}

The same argument that states that $\mathbf{U^*XV}$ has to be
diagonal, can also be applied when the constraint is given by the
nuclear norm. Define $\tilde{\mathbf{\Sigma}}=\mathbf{U^*XV}$. We
wish to minimize $\Vert
\mathbf{\tilde{\mathbf{\Sigma}}}-\mathbf{\Sigma}
\Vert_F=\sum_i{(\tilde{\sigma_i}-\sigma_i)^2}$ s.t. $\Vert
\mathbf{X} \Vert_*=\Vert \mathbf{\tilde{\Sigma}} \Vert_*=\sum_i{|
\tilde{\sigma_i}} | \le \lambda , i=1, \ldots k, k \le n$. Note that
$\tilde{\sigma_i}$ has to be nonnegative otherwise it will increase
the Frobenius norm but will not change the nuclear norm. Hence, the
problem can now be formulated as:
\begin{equation}
\label{eq_nuclear}
\begin{array}{l}
\mbox{minimize   } \sum_i{(\tilde{\sigma_i}-\sigma_i)^2} \\
\mbox{subject to } \sum_i{\tilde{\sigma_i}} \le \lambda, \\
\mbox{~~~~~~~~~~~~~} \tilde{\sigma_i}\ge 0 .
 \end{array}
\end{equation}

This is a standard convex optimization problem that can be solved by
methods such as semidefinite programming \cite{boyd}. The exact same 
can be done to the Ky-Fan norm.

\section{Approximation of certain entries}
Suppose we wish to approximate only certain entries of the matrix, under different constraints, i.e. we 
are interested in solving Eq.~\ref{eq1}, given that the solution of Eq.~\ref{eq2} is known and given by 
$\mathcal{D}\mathbf{M}$, where $\mathcal{D}$ is the solution operator. For example, if the constraint is $\text{rank}(\mathbf{X})\le k$
$\mathcal{D}\mathbf{X}$ is the truncated SVD of $\mathbf{X}$ containing the first $k$ singular values. Note that $\mathcal{D}$ is not necessarily convex.
We examine the following iterative algorithm:
\begin{equation}
\label{projgrad}
\mathbf{X}_{n+1}=\mathcal{D}(\mathbf{X}_n-\mathcal{P}(\mathbf{X}_n-\mathbf{M})).
\end{equation}
Eq.~\ref{projgrad} can be considered as a projected gradient algorithm with unit step size, where the projection is
given by $\mathbf{D}$. 
\begin{theorem} [Local Convergence]: Let $\epsilon(\mathbf{X}_n)=\Vert \mathcal{P}\mathbf{X}_n-\mathcal{P}\mathbf{M}\Vert_F$ be the error 
at the $n$th iteration, then $\epsilon(\mathbf{X}_n)$ is monotonically decreasing, and because it is bounded the algorithm converges.
\label{th1}
\end{theorem}
The proof for Theorem~\ref{th1} is given in \cite{gil}. Theorem~\ref{th1} does not say anything about convergence to the global
solution. However, when the projection $\mathcal{D}$ is convex and self adjoint ($\mathcal{D=D^*}$) and the algorithm is modified to have adaptive step size, that is:
\begin{equation}
\label{projgrad_adap}
\mathbf{X}_{n+1}=\mathcal{D}(\mathbf{X}_n-\mu_n\mathcal{P}(\mathbf{X}_n-\mathbf{M})),
\end{equation}
and $\mu_n$ is computed by Armijo rule in a greedy form, minimizing the error in every iteration:

\begin{equation}
\label{armijo}
\begin{array}{l}
l[n]=\mbox{argmin}_{j\in \mathcal{Z}_{\ge 0}}: f(\mathbf{X}_{n,j}) \\
\le f(\mathbf{X}_n)- \sigma \mbox{trace}(\nabla f(\mathbf{X}_n)^T (\mathbf{X}_n-\mathbf{Z}_{n,j})) \\
\mathbf{Z}_{n,j}=\mathcal{D}(\mathbf{X}_n-\tilde{\mu}2^{-j}\nabla f(\mathbf{X}^n)) \\
\mu_n=\tilde{\mu}2^{-l[n]},
\end{array}
\end{equation}
where $f(X)=\frac{1}{2}\Vert \mathcal{P}\mathbf{X}-\mathcal{P}\mathbf{M} \Vert_F^2$, $\tilde{\mu}>0$ and $\sigma \in (0,1)$, Then the algorithm is guarantee to achieve the global solution \cite{iusem}. This approach has two major problems:
\begin{itemize}
\item
For the cases of interest, the operators for truncating the nuclear and spectral norm, are not self-adjoint ($\mathcal{D \neq D^*}$)
\item
This approach requires applying the Armijo rule in every iteration. This means several applications of the operator $\mathcal{D}$ in each iteration which is usually computationally expensive.
\end{itemize}

As for the first point, requiring the projection $\mathcal{D}$ to be self-adjoint can be slightly more than needed for the global convergence proof in \cite{iusem}. This requirement is needed in order to satisfy $\langle X-Y, \mathcal{D}X-X \rangle \geq 0$ for $Y=\mathcal{D}Y$, which always holds when $\mathcal{D=D^*}$, but also when $\mathcal{D}$ is as we defined in Lemma \ref{spectral} and Eq. \ref{eq_nuclear}.

\begin{theorem}
Let $\mathcal{D}$ be the following projection (defined as in Lemma \ref{spectral}): Given the SVD of X is $\mathbf{X=USV}^*$, we define $\mathcal{D}_\lambda X=U\tilde{S}V^*$ where $\tilde{s}_i=$min$(s_i,\lambda)$. Then, for every matrices $\mathbf{X}$ and $\mathbf{Y}$ such that $Y=\mathcal{D}\mathbf{Y}$, $\langle \mathbf{X}-\mathbf{Y}, \mathcal{D}\mathbf{X}-\mathbf{X} \rangle \geq 0$
\end{theorem}
\begin{IEEEproof}
The condition $\langle \mathbf{X}-\mathbf{Y}, \mathcal{D}\mathbf{X}-\mathbf{X} \rangle \geq 0$ can be reformulated as
\begin{equation}
\label{bb}
\langle \mathbf{X}, \mathbf{X}-\mathcal{D}\mathbf{X} \rangle \geq \langle \mathbf{Y}, \mathbf{X}-\mathcal{D}\mathbf{X} \rangle,
\end{equation}
where $\Vert Y \Vert_2 \le \lambda$.

First, note that according to the Von Neumann's trace inequality, the value of the right hand side is maximal when $\mathbf{Y}$ and $\mathbf{X}-\mathcal{D}\mathbf{X}$ have the singular vectors.
Hence, we define: $\mathbf{X}=\mathbf{US}_X\mathbf{V}^*$, $\mathbf{Y}=\mathbf{U}\tilde{\mathbf{S}}_Y\mathbf{V}^*$ and $\mathcal{D}\mathbf{X}=\mathbf{U}\tilde{\mathbf{S}}_X\mathbf{V}^*$. The tilde is for indicating that the singular values of $\tilde{\mathbf{S}}$ are smaller or equal to $\lambda$. 

We start by evaluating the left side of Eq.\ref{bb}:
\begin{equation}
\langle \mathbf{X}, \mathbf{X}-\mathcal{D}\mathbf{X} \rangle=\text{trace}[\mathbf{S}_X(\mathbf{S}_X-\tilde{\mathbf{S}_X})]=\sum_i s_{x_i}(s_{x_i}-\tilde{s}_{x_i}).
\end{equation}
Now, for $s_{x_i} \le \lambda$ we get $(s_{x_i}-\tilde{s}_{x_i})=0$. Hence, only when $s_{x_i} > \lambda$ the sum grows and the expression can be rewritten as: $\langle \mathbf{X}, \mathbf{X}-\mathcal{D}\mathbf{X} \rangle=\sum_{s_{x_i}>\lambda} s_{x_i}(s_{x_i}-\tilde{s}_{x_i})$

We now observe the right side of Eq. \ref{bb}:
\begin{equation}
\langle \mathbf{Y}, \mathbf{X}-\mathcal{D}\mathbf{X} \rangle=\text{trace}[\tilde{\mathbf{S}}_Y(\mathbf{S}_X-\tilde{\mathbf{S}}_X)]=\sum_i \tilde{s}_{y_i}(s_{x_i}-\tilde{s}_{x_i}).
\end{equation}
Again, the elements that contribute to the sum are those for which $s_{x_i} > \lambda$. Hence, on the right side we obtained: $\langle \mathbf{Y}, \mathbf{X}-\mathcal{D}\mathbf{X} \rangle=\sum_{s_{x_i}>\lambda} \tilde{s}_{y_i}(s_{x_i}-\tilde{s}_{x_i})$.

Both expressions can be thought of as a sum of the positive elements $(s_{x_i}-\tilde{s}_{x_i})$ with different coefficients. Both series have the same length ($s_{x_i} > \lambda$) but the coefficient on the left side is $s_{x_i}$ for $i$'s that give $s_{x_i}>\lambda$ and the right hand series coefficients are by definition (since $\Vert \mathbf{Y} \Vert_2 \le \lambda$)  smaller than $\lambda$. Therefore, the sum of the left side is bigger than the sum of the right side. This completes the proof.
\end{IEEEproof}
This means that for the spectral norm, the algorithm converges to the global solution. The exact same proof can be done for the nuclear norm and Ky-Fan norm as well, showing the algorithm converges to global solution.

\begin{theorem}[Optimal step size]
For the matrix approximation problem (Eq.~\ref{eq1}) with convex $\mathcal{D}$, the optimal step size is given by $\mu_n=1$.
\label{th2}
\end{theorem}
The proof of Theorem~\ref{th2} is given in \cite{gil}. Note that this holds for any case of projected gradient involving 
orthogonal axes. Theorem~\ref{th2} states that in our case, when having a convex constraint and projection, then Eq.~\ref{projgrad} converges to the global solution. This means, that now we can solve a variety of matrix approximation problem with reasonable computation rate. Note, that we have shown that in some cases, global solution is achieved even when the projection is not self-adjoint (orthogonal).
The next section shows, how this very simple algorithm, can be applied to matrix completion problems as well.

\section{Matrix Completion}
Matrix completion is an important problem that has been investigated
extensively. The matrix completion problem differs from the
matrix approximation problem by the fact that the known entries must
remain fixed while changing their role from the objective function
to be minimized to the constraint part. A well investigated  matrix
completion problem appears in the introduction as the rank
minimization problem. Because rank minimization is not convex and
NP-hard, it is usually relaxed for the nuclear norm minimization.
Since for the convex case, we have seen that Eq.~\ref{projgrad} converges to the
global solution, matrix completion can be achieved simply by using binary
search. The advantage of this approach over other different
approaches, which minimize the nuclear norm for example, is that it
is general and can be applied to other problems that were not
addressed such as minimizing the spectral norm. Moreover, some
algorithms such as the Singular Value Thresholding (SVT)
\cite{candes} require additional parameters $\tau$ and $\delta$ that
affect the convergence and the final result, where in this approach
no external parameters are required (except for tolerance for determining convergence).

\begin{algorithm}[!ht]
\caption{Matrix Completion using  Nuclear Norm / Spectral Norm Minimization}
\label{compalgo}
\textbf{Input:}  $\mathbf{M}$ - matrix to
complete, $\mathcal{P}$ - projection operator that specifies the
important entries, \\ $tol$ - admissible approximation error, $\lambda_{tol}$ - admissible constraint accuracy  \\ \protect
\textbf{Output:} $\mathbf{X}$ - Completed matrix \protect
\begin{algorithmic}[1]
\STATE $\mathbf{M} \leftarrow \mathcal{P}\mathbf{M}$ 
\STATE $\lambda_{min}\leftarrow 0$
\STATE $\lambda_{max}\leftarrow \Vert \mathbf{M} \Vert_*$ (or $\Vert \mathbf{M} \Vert_2$ for the spectral norm)
\STATE $\lambda \leftarrow 0$
\REPEAT
\STATE $\lambda_{prev} \leftarrow \lambda$
\STATE $\lambda \leftarrow (\lambda_{min}+\lambda_{max})/2$
\STATE $\mathbf{X} \leftarrow$ Approximate $\mathcal{P}\mathbf{M}$ s.t. $\Vert \mathbf{X} \Vert_* \le \lambda$ (or $\Vert \mathbf{X} \Vert_2 \le \lambda$ for the spectral norm case)
\STATE $error \leftarrow \Vert \mathcal{P}\mathbf{X}-\mathcal{P}\mathbf{M}\Vert_F$
\IF{$error > tol$} 
\STATE $\lambda_{min}\leftarrow \lambda$
\ELSE 
\STATE $\lambda_{max}\leftarrow \lambda$
\ENDIF

\UNTIL $error<tol$ and $|\lambda-\lambda_{prev}| < \lambda_{tol}$
\RETURN $\mathbf{X}$
\end{algorithmic}
\end{algorithm}

This approach is detailed in Algorithm~\ref{compalgo}, which is robust and does not require any tuning, other than tolerance threshold for determining convergence.
Algorithm~\ref{compalgo} can be used for a matrix completion under a variety of constraints.

\begin{figure}
	\centering
		\includegraphics[width=0.45\textwidth]{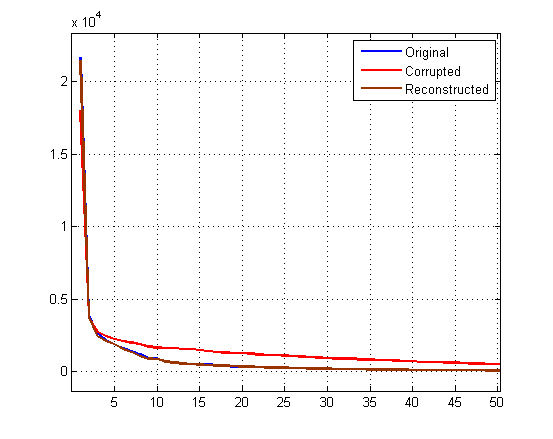}
	\label{fig:svd_comp}
	\caption{Singular values comparison between the different images.}
\end{figure}

\begin{figure}
	\centering
		\includegraphics[width=0.24\textwidth]{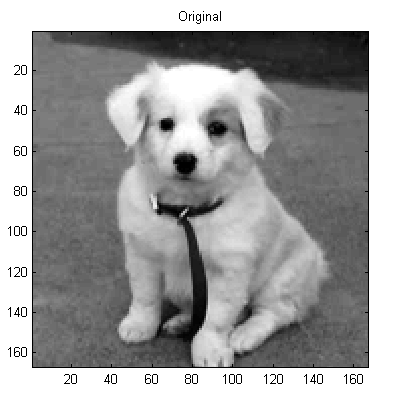}
		\includegraphics[width=0.24\textwidth]{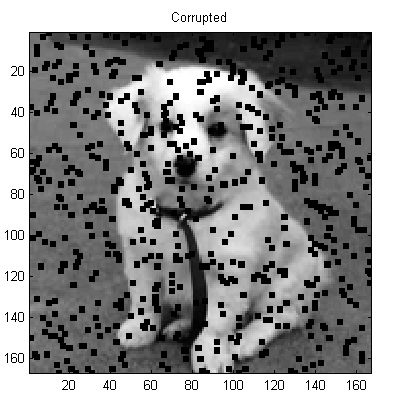}
		\includegraphics[width=0.24\textwidth]{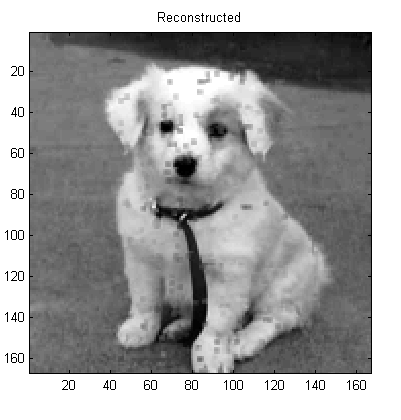}
	\label{fig:corrupted_dog}
	\caption{Corrupted dog image and the reconstructed image.}
\end{figure}

Fig. \ref{fig:corrupted_dog} shows Algorithm~\ref{compalgo} results over a corrupted image. In the corrupted image, squares of size $3 \times 3$ were randomly removed from the image, destroying $18\%$ of it. The reconstruction is more difficult, since the damage is in squares and not just irregular points. The original image nuclear norm is $51,625$, the corrupted nuclear norm is $96,500$ and the norm of the completed matrix is $50,418$. Minimizing nuclear norm for image reconstructing is a well known method, as images usually have a low numerical rank as the singular values decay very fast. It can be seen in Fig. \ref{fig:svd_comp} that the singular values of the reconstructed image, are almost identical to the original.

\section*{Acknowledgment}
\noindent This research was partially supported by the Israel Science
Foundation (Grant No. 1041/10) and by the Israeli Ministry of Science \& Technology
3-9096.


\begin{thebibliography}{99}
\bibitem{sta} T.A. Louis, \emph{Finding the observed information matrix when using the EM algorithm}, Journal of the Royal Statistical Society, Series B. (Methodological), Vol. 44, No. 2, pp. 226-233, 1982.
\bibitem{ml} T. Hastie, R. Tibshirani, G. Sherlock, M. Eisen, P. Brown, and D. Botstein, \emph{Imputing missing data for gene expression arrays}. Technical report; Division of Biostatistics, Stanford University, 1999.
\bibitem{SJ} N. Srebro and T. Jaakkola, \emph{Weighted low-rank approximations}, Preceeding of the 20th International Conference on Machine Learning (ICML-2003), Washington DC, 2003.
\bibitem{mikielad} J. Mairal, M. Elad, G. Sapiro, \emph{Sparse representation for color image restoration}, IEEE Transactions on Image Processing, Vol. 17, No. 1, pp.53-69, 2008.
\bibitem{fazel} M. Fazel, \emph{Matrix Rank Minimization with Applications. PhD thesis}, Stanford University, 2002.
\bibitem{candes} J.F. Cai, E.J. Candes and Z. Shen, \emph{Singular Value Thresholding Algorithm for Matrix Completion}, SIAM Journal on Optimization,  Vol. 20, No. 4, pp. 1956-1982, 2010.
\bibitem{Nan} Feng Nan, \emph{Low Rank Matrix Completion}, Master thesis, Massachusetts Institute of Technology, 2009.
\bibitem{peter} P. H. Schonemann, \emph{A generalized solution of the orthogonal procrustes problem}, Psychometrika, Vol. 31, No. 1, pp. 1-10, 1966.
\bibitem{everson} R. Everson, \emph{Orthogonal but not orthonormal Procrustes problem}, 1997.
\bibitem{bhatia} R. Bhatia, \emph{Matrix Analysis}, Graduate Texts in Mathematics, Springer 1996.
\bibitem{boyd} S. Boyd, L. Vandenberghe, \emph{Convex Optimization}, Cambridge University Press, 2004.
\bibitem{goldberg}, I.C. Gohberg, M.G. Krein, \emph{Introduction to the theory of linear and selfadjoint operators}, Translations of Mathematical Monographs, Vol. 18, pp. 94-95, 1969.
\bibitem{iusem} A.N. Iusem, \emph{On the convergence properties of the projected gradient method for convex optimization}, Computational and Applied Mathematics, Vol. 22, No. 1, pp. 37-52, 2003.
\bibitem{gil} G. Shabat, A. Averbuch, \emph{Interest Zone Matrix Approximation}, Electronic Journal of Linear Algebra, Vol. 23, pp. 678-702, 2012.
\end{thebibliography}
\end{document}